\documentclass{elsart}
\usepackage{amssymb,amsfonts}

\newcommand{\dom}{\mathop{\rm dom}\nolimits}
\newcommand{\trace}{\mathop{\rm Tr}\nolimits}

\newcommand{\twomat}[4]{\left(\begin{array}{cc}#1&#2\\#3&#4\end{array}\right)}

\newcommand{\cA}{{\mathcal A}} 
\newcommand{\cB}{{\mathcal B}} 

\DeclareRobustCommand\openone{\leavevmode\hbox{\small1\normalsize\kern-.33em1}}
\newcommand{\id}{\mathbf{I}}

\newcommand{\be}{\begin{equation}}
\newcommand{\ee}{\end{equation}}
\newcommand{\bea}{\begin{eqnarray}}
\newcommand{\eea}{\end{eqnarray}}
\newcommand{\beas}{\begin{eqnarray*}}
\newcommand{\eeas}{\end{eqnarray*}}
\newtheorem{definition}{Definition}
\newtheorem{theorem}{Theorem}
\newtheorem{lemma}{Lemma}

\newtheorem{proposition}{Proposition}
\newcount\minute
\newcount\hour
\def\currenttime{%
    \minute\time
    \hour\minute
    \divide\hour60
    \the\hour:\multiply\hour60\advance\minute-\hour\the\minute}
\begin{document}
\begin{frontmatter}
\title{An Araki-Lieb-Thirring inequality for geometrically concave and geometrically convex functions}
\author{Koenraad M.R. Audenaert}
\address{
Department of Mathematics, Royal Holloway, University of London,\\
Egham TW20 0EX, United Kingdom}
\ead{koenraad.audenaert@rhul.ac.uk}
\date{\today, \currenttime}
\begin{abstract}
For positive definite matrices $A$ and $B$, the Araki-Lieb-Thirring inequality amounts to an
eigenvalue log-submajorisation relation for fractional powers
$$
\lambda(A^t B^t) \prec_{w(\log)} \lambda^t(AB), \quad 0<t\le 1,
$$
while for $t\ge1$, the reversed inequality holds.
In this paper I generalise this inequality, replacing the fractional powers $x^t$ by a larger class of functions.
Namely, a continuous, non-negative, geometrically concave function
$f$ with domain $\dom(f)=[0,x_0)$ for some positive $x_0$ (possibly infinity) satisfies
$$
\lambda(f(A) f(B)) \prec_{w(\log)} f^2(\lambda^{1/2}(AB)),
$$
for all positive semidefinite $A$ and $B$ with spectrum in $\dom(f)$,
if and only if $0\le xf'(x)\le f(x)$ for all $x\in\dom(f)$.
The reversed inequality holds for continuous, non-negative, geometrically convex functions if and only if they satisfy
$xf'(x)\ge f(x)$ for all $x\in\dom(f)$.
As an application I derive a complementary inequality to the Golden-Thompson inequality.
\end{abstract}

\end{frontmatter}
\section{Introduction\label{sec1}}
The Araki-Lieb-Thirring (ALT) inequality \cite{araki,lieb} states
that for $0<t\le 1$ and positive definite matrices $A$ and $B$,
the eigenvalues of $A^tB^t$ are log-submajorised by the eigenvalues of $(AB)^t$.
For $n\times n$ matrices $X$ and $Y$ with positive spectrum,
the log-submajorisation relation $\lambda(X)\prec_{w(\log)}\lambda(Y)$ means
that for all $k=1,\ldots,n$, the following holds:
$$
\prod_{j=1}^k\lambda_j(X) \le \prod_{j=1}^k\lambda_j(Y).
$$
This is equivalent to weak majorisation of the logarithms of the spectra
$$
\log\lambda(X) \prec_w \log\lambda(Y).
$$
Here and elsewhere, I adhere to the convention to sort eigenvalues in non-increasing order;
that is, $\lambda_1(X)\ge\lambda_2(X)\ge\ldots\ge\lambda_n(X)$.


With this notation, the ALT inequality can be written as
\be
\lambda(A^tB^t) \prec_{w(\log)} \lambda^t(AB),\quad 0<t\le 1. \label{eq:ALT}
\ee
For positive scalars $a$ and $b$, (\ref{eq:ALT}) reduces to the equality $a^tb^t = (ab)^t$.

One can ask whether similar inequalities hold for other functions than the fractional powers $x^t$.
One possibility is to consider functions that satisfy
$$
\lambda(f(A)f(B)) \prec_{w(\log)} f(\lambda(AB)).
$$
As the scalar case reduces to $f(a)f(b)\le f(ab)$ these functions must be super-multiplicative.
Another possibility, and the one pursued here, is to consider functions satisfying
\be
\lambda(f(A)f(B)) \prec_{w(\log)} f^2\left(\sqrt{\lambda(AB)}\right).\label{eq:falt}
\ee
Here, the scalar case reduces to $f(a)f(b)\le f^2(\sqrt{ab})$, for all
$a,b>0$.
Functions satisfying this requirement are called \textit{geometrically concave} (see Definition \ref{def:geo} below).
In this paper I completely characterise the class
of geometrically concave functions that satisfy (\ref{eq:falt}) for all positive definite matrices
$A$ and $B$.

Likewise, as inequality (\ref{eq:ALT}) holds in the reversed sense for $t\ge 1$, one may ask for which functions $f$
the reversed inequality holds for all positive definite matrices
$A$ and $B$:
\be
 f^2\left(\sqrt{\lambda(AB)}\right) \prec_{w(\log)} \lambda(f(A)f(B)).\label{eq:ralt}
\ee
Here the scalar case restricts the class of functions to those satisfying the relation
$f^2(\sqrt{ab})\le f(a)f(b)$. Such functions are called \textit{geometrically convex}.
I also completely characterise the class
of geometrically convex functions that satisfy (\ref{eq:ralt}) for all positive definite matrices
$A$ and $B$.

The concepts of geometric concavity and geometric convexity were first studied by Montel \cite{montel}
and have recently received attention from the matrix community \cite{ka,bourin}.
\begin{definition}\label{def:geo}
Let $I$ be the interval $I=[0,x_0)$, with $x_0>0$ (possibly infinite).
A function $f:I\to[0,\infty)$ is \emph{geometrically concave} if for all $x,y\in I$,
$\sqrt{f(x)f(y)} \le f(\sqrt{xy})$.
It is \emph{geometrically convex} if for all $x,y\in I$,
$\sqrt{f(x)f(y)} \ge f(\sqrt{xy})$.
\end{definition}
Equivalently, a function $f(x)$ is geometrically concave (convex)
if and only if the \textit{associated function} $F(y):=\log(f(e^y))$ is concave (convex).





The main results of this paper are summarised in the next section,
the proofs of the main theorems (Theorems \ref{th:FALT} and \ref{cor:RALT}) are given in Section \ref{sec:proof},
and the paper concludes with a brief application in Section \ref{sec3}.
\section{Main Results}
To state the main results of this paper most succinctly, let us define two classes of functions.
\begin{definition}
A continuous non-negative function $f$ with domain an interval $I=[0,x_0)$ is in class $\cA$ if and only if
it is geometrically concave and its derivative $f'$ satisfies $0\le xf'(x)\le f(x)$ for all $x\in I$ where
the derivative exists.
\end{definition}
\begin{definition}
A continuous non-negative function $f$ with domain an interval $I=[0,x_0)$ is in class $\cB$ if and only if
it is geometrically convex and its derivative $f'$ satisfies $xf'(x)\ge f(x)$ for all $x\in I$ where
the derivative exists.
\end{definition}
In terms of the associated function $F(y)=\log(f(\exp y))$, $f\in\cA$ if and only if $F(y)$ is concave
and $0\le F'(y)\le 1$ for all $y$ where $F$ is differentiable, and $f\in\cB$ if and only if
$F(y)$ is convex and $1\le F'(y)$ for all $y$ where $F$ is differentiable.

There is a simple one-to-one relationship between these two classes; 
essentially $f$ is in class $\cA$ if and only if its inverse function $f^{-1}$ is in class $\cB$.
However, some care must be taken as $\cA$ contains the constant functions and also those functions that are constant
on some interval.
\begin{proposition}
A function $f$ that is non-constant on the interval $[0,x_1]\subseteq[0,x_0)$ is in class $\cA$ if and only if 
the inverse of the restriction of $f$ to $[0,x_1]$ is in class $\cB$.
\end{proposition}
\textit{Proof.}
It is clear that a concave monotonous function $f$ is always invertible over the entire interval where it is not constant.
We will henceforth identify the inverse of $f$ with the inverse of the restriction of $f$ on that interval.

If $F(y)$ is the associated function of $f(x)$
then the associated function of $f^{-1}$ is the inverse function of $F$, $F^{-1}$.
Now $f$ is in class $\cA$ if and only if $F$ is concave, monotonous and $F'\le 1$.
This implies that the inverse function $G=F^{-1}$ is convex and satisfies $G'\ge1$,
which in turn implies that $G$ is the associated function of a function $g$
in class $\cB$.
This shows that $f\in\cA$ implies $f^{-1}\in\cB$.

A similar argument reveals that the converse statement holds as well.
\qed

The main result of this paper is the following theorem:
\begin{theorem}\label{th:FALT}
Let $f$ be a continuous non-negative function with domain an interval $I=[0,x_0)$, $x_0>0$ (possibly infinite),
then
\be
\lambda(f(A)f(B))\prec_{w(\log)} f^2\left(\sqrt{\lambda(AB)}\right)\label{eq:FALT}
\ee
holds for all positive definite matrices $A$ and $B$ with spectrum in $I$
if and only if $f$ is in class $\cA$.
\end{theorem}
That the right-hand side of (\ref{eq:FALT}) is well-defined follows from the following lemma:
\begin{lemma}\label{lem:a1}
If $A$ and $B$ are positive semidefinite matrices with eigenvalues in the interval $I=[a,b]$, $0\le a<b$,
the positive square roots of the eigenvalues of $AB$ are in $I$ as well.
\end{lemma}
\textit{Proof.}
We have $a\le A,B\le b$, which implies
$$
a^2\le aA \le A^{1/2}BA^{1/2} \le bA\le b^2.
$$
Hence, $a^2\le \lambda_i(A^{1/2}BA^{1/2})\le b^2$, so that
$a\le \lambda_i^{1/2}(AB)\le b$.
\qed

A simple consequence of Theorem \ref{th:FALT} is that the reversed inequality holds if and only if $f$ is in class $\cB$.
\begin{theorem}\label{cor:RALT}
Let $g$ be a continuous non-negative function with domain an interval $I=[0,x_0)$, $x_0>0$ (possibly infinite),
then
\be
 g^2\left(\sqrt{\lambda(XY)}\right) \prec_{w(\log)} \lambda(g(X)g(Y)).\label{eq:RALT}
\ee
holds for all positive definite matrices $X$ and $Y$ with spectrum in $I$
if and only if $g$ is in class $\cB$.
\end{theorem}

\section{Proofs\label{sec:proof}}
We now turn to the proofs of Theorems \ref{th:FALT} and \ref{cor:RALT}.

\subsection{Proof of necessity}
To show necessity of the conditions $f\in\cA$ ($f\in\cB$)
I consider two special $2\times 2$ matrices with eigenvalues $a$ and $b$, $0\le b<a$, 
such that $a\in\dom(f)$ and $f$ is differentiable in $a$: 
$$
A=\twomat{a}{0}{0}{b}, \quad B=U\twomat{a}{0}{0}{b}U^*,\mbox{ with }U=\twomat{\cos\theta}{-\sin\theta}{\sin\theta}{\cos\theta}.
$$
We will consider values of $\theta$ close to 0.
The largest eigenvalue of $AB$ can be calculated in a straight-forward fashion.
The quantity $f^2(\sqrt{\lambda_1(AB)})$ can then be expanded in a power series of the variable $\theta$. To second order this
yields
$$
f^2(\sqrt{\lambda_1(AB)}) = f(a)^2 -\frac{a(a-b)f(a)f'(a)}{a+b}\theta^2+O(\theta^4).
$$
In a similar way we also get
$$
\lambda_1(f(A)f(B)) = f(a)^2-\frac{f(a)^2(f(a)-f(b))}{f(a)+f(b)}\theta^2+O(\theta^4).
$$
Hence, to satisfy inequality (\ref{eq:FALT}), the following must be satisfied for all $0\le b<a\in\dom(f)$:
$$
\frac{a(a-b)f(a)f'(a)}{a+b} \le \frac{f(a)^2(f(a)-f(b))}{f(a)+f(b)}.
$$
In particular, take $b=0$. As $f$ has to be geometrically concave, $f(0)=0$.
The condition then becomes
$$
af(a)f'(a) \le f(a)^2, \forall a\in\dom(f),
$$
which reduces to the defining condition for $f\in\cA$.

Necessity of the condition $xf'(x)\ge f(x)$ in Corollary \ref{cor:RALT}
also follows immediately from this special pair of matrices.

Note that in the preceding proof we see why the domain of $f$ should include the point $x=0$.

\subsection{Proof of sufficiency for Theorem \ref{th:FALT}}
Now I turn to proving sufficiency of the condition $f\in\cA$.
The main step consists in showing that the set of
functions $f$ for which the inequality (\ref{eq:FALT}) holds is `geometrically convex'; that is, the set
of associated functions $F$ for these $f$ is convex. To show this a number of preliminary propositions are needed.

\begin{lemma}\label{lem:u1}
Let $R_1$, $R_2$, $S_1$ and $S_2$ be positive semidefinite matrices such that $R_1$ commutes with $R_2$
and $S_1$ with $S_2$.
Let $R=R_1^{1/2}R_2^{1/2}$ and $S=S_1^{1/2}S_2^{1/2}$.
Then
\be
\lambda_1((R^{1/2}SR^{1/2})^2) \le \lambda_1(R_1^{1/2}S_1R_1^{1/2}\;R_2^{1/2}S_2R_2^{1/2}).\label{eq:u1}
\ee
\end{lemma}
\textit{Proof.}
We will prove this by showing that the equality
$$\lambda_1(R_1^{1/2}S_1R_1^{1/2}\;R_2^{1/2}S_2R_2^{1/2})= a$$
implies
the inequality $\lambda_1((R^{1/2}SR^{1/2})^2)\le a$.

W.l.o.g.\ we can assume that the matrices $R_1$ and $S_1$ are invertible;
then the equality indeed leads to the following sequence of implications:
\beas
&&\lambda_1(R_1^{1/2}S_1R_1^{1/2}\;R_2^{1/2}S_2R_2^{1/2})= a\\
&\Longrightarrow&(R_1^{1/2}S_1R_1^{1/2})^{1/2}\;R_2^{1/2}S_2R_2^{1/2}\;(R_1^{1/2}S_1R_1^{1/2})^{1/2} \le a \\
&\Longrightarrow& R_2^{1/2}S_2R_2^{1/2} \le a\;R_1^{-1/2}S_1^{-1}R_1^{-1/2} \\
&\Longrightarrow& S_1^{1/2}R_1^{1/2}\;R_2^{1/2}S_2R_2^{1/2}\;R_1^{1/2}S_1^{1/2} =
S_1^{1/2}RS_2RS_1^{1/2} \le a \\
&\Longrightarrow& \sigma_1(S_2^{1/2}RS_1^{1/2}) \le \sqrt{a} \\
&\Longrightarrow& |\lambda_1(S_2^{1/2}RS_1^{1/2})| \le \sqrt{a}.
\eeas
The last implication is the simplest case of Weyl's majorant theorem.

Now note that $S_2^{1/2}RS_1^{1/2}$ and $S^{1/2}RS^{1/2}=S_1^{1/4}S_2^{1/4}RS_1^{1/4}S_2^{1/4}$
have the same non-zero eigenvalues.
Hence, $\lambda_1(S^{1/2}RS^{1/2})=\lambda_1(R^{1/2}SR^{1/2}) \le \sqrt{a}$,
and the inequality $\lambda_1((R^{1/2}SR^{1/2})^2)\le a$ follows.
\qed

\begin{proposition}\label{cor:u2}
Under the conditions of Lemma \ref{lem:u1},
\be
\prod_{i=1}^k\lambda_i(R_1^{1/2}R_2^{1/2}\;S_1^{1/2}S_2^{1/2})
\le \prod_{i=1}^k\lambda_i^{1/2}(R_1S_1)\;\lambda_i^{1/2}(S_2R_2).\label{eq:u2}
\ee
\end{proposition}
\textit{Proof.}
Since each side of (\ref{eq:u1}) is the largest eigenvalue of a product of powers of matrices, we can use
the well-known Weyl trick of replacing every matrix by its antisymmetric tensor power
to boost the inequality to the log-submajorisation relation
$$
\prod_{i=1}^k\lambda_i((R^{1/2}SR^{1/2})^2) \le \prod_{i=1}^k\lambda_i(R_1^{1/2}S_1R_1^{1/2}\;R_2^{1/2}S_2R_2^{1/2}).
$$
Combining this with Lidskii's inequality $\prod_{i=1}^k\lambda_i(AB)\le \prod_{i=1}^k\lambda_i(A)\lambda_i(B)$
(\cite{bhatia}, Corollary III.4.6),
valid for positive definite $A$ and $B$, and then taking square roots yields inequality (\ref{eq:u2}).
\qed

Inequality (\ref{eq:u2}) can be interpreted as midpoint geometric convexity of the function
$$
p\mapsto f_k(p)=\prod_{i=1}^k\lambda_i(R_1^{p}R_2^{1-p}\;S_1^{p}S_2^{1-p});
$$
that is, $f_k(1/2)\le \sqrt{f_k(1)f_k(0)}$.
We now use a standard argument (see e.g.\ the proof of Lemma IX.6.2 in \cite{bhatia})
to show that this actually implies geometric convexity in full generality,
i.e.\ $f_k(p)\le f_k(1)^p f_k(0)^{1-p}$ for all $p\in[0,1]$.
\begin{proposition}\label{cor:u3}
Under the conditions of Lemma \ref{lem:u1}, and for all $p\in[0,1]$,
\be
\prod_{i=1}^k\lambda_i(R_1^{p}R_2^{1-p}\;S_1^{p}S_2^{1-p})
\le \prod_{i=1}^k\lambda_i^{p}(R_1S_1)\;\lambda_i^{1-p}(S_2R_2).\label{eq:u3}
\ee
\end{proposition}
\textit{Proof.}
By Proposition \ref{cor:u2} the inequality holds for $p=1/2$. It trivially holds for $p=0$ and $p=1$.

Let $s,t\in[0,1]$ be given.
Applying Proposition \ref{cor:u2} with the matrices
$R_1$, $S_1$, $R_2$ and $S_2$ replaced by
$R_1^t R_2^{1-t}$, $S_1^t S_2^{1-t}$,
$R_1^s R_2^{1-s}$ and $S_1^s S_2^{1-s}$,
respectively, yields the inequality
\beas
\lefteqn{\prod_{i=1}^k\lambda_i(R_1^{(s+t)/2}R_2^{1-(s+t)/2}\;S_1^{(s+t)/2}S_2^{1-(s+t)/2})} \\
&\le& \prod_{i=1}^k\lambda_i^{1/2}(R_1^tR_2^{1-t} S_1^tS_2^{1-t})\;\lambda_i^{1/2}(R_1^sR_2^{1-s} S_1^sS_2^{1-s}).
\eeas
Now assume that the inequality (\ref{eq:u3}) holds for the values $p=s$ and $p=t$.
Thus
\beas
\lefteqn{\prod_{i=1}^k\lambda_i^{1/2}(R_1^tR_2^{1-t} S_1^tS_2^{1-t})\;\lambda_i^{1/2}(R_1^sR_2^{1-s} S_1^sS_2^{1-s})} \\
&\le&
\prod_{i=1}^k\lambda_i^{t/2}(R_1S_1)\;\lambda_i^{(1-t)/2}(S_2R_2)\;\lambda_i^{s/2}(R_1S_1)\;\lambda_i^{(1-s)/2}(S_2R_2)\\
&=& \prod_{i=1}^k\lambda_i^{(s+t)/2}(R_1S_1)\;\lambda_i^{1-(s+t)/2}(S_2R_2).
\eeas
In other words, the assumption that (\ref{eq:u3}) holds for the values $p=s$ and $p=t$ implies that it also holds
for their midpoint $p=(s+t)/2$.

Using induction this shows that (\ref{eq:u3}) holds for all dyadic rational values of $p$
(i.e.\ rationals of the form $k/2^n$, with $k$ and $n$ integers such that $k\le 2^n$). Invoking continuity
and the fact that the dyadic rationals are dense in $[0,1]$,
this finally implies that (\ref{eq:u3}) holds for all real values of $p$ in $[0,1]$.
\qed

We are now ready to prove our first intermediate result: convexity of the set of associated functions $F$ for which
the inequality (\ref{eq:FALT}) holds.
\begin{proposition}\label{prop:convex}
Let $f_1(x)$ and $f_2(x)$ be two continuous, non-negative functions with domain an interval $I$ of the non-negative reals, and
for which (\ref{eq:FALT}) holds for all positive semidefinite $A$ and $B$
with spectrum in $I$.
Let $p\in[0,1]$ and let $f(x)=f_1^p(x)f_2^{1-p}(x)$. Then (\ref{eq:FALT}) holds for $f$ too.
\end{proposition}
\textit{Proof.}
Let us fix the matrices $A$ and $B$ and let $R_i=f_i(A)$ and $S_i=f_i(B)$, $i=1,2$.
These matrices $R_i$ and $S_i$ clearly satisfy the conditions of Proposition \ref{cor:u3} (positivity and commutativity).
Hence
\beas
\prod_{i=1}^k\lambda_i(f(A)\;f(B))
&=&
\prod_{i=1}^k\lambda_i(f_1^{p}(A)f_2^{1-p}(A)\;f_1^{p}(B)f_2^{1-p}(B))\\
&\le& \prod_{i=1}^k\lambda_i^{p}(f_1(A)f_1(B))\;\lambda_i^{1-p}(f_2(A)f_2(B)).
\eeas
By the assumption that $f_1$ and $f_2$ satisfy inequality (\ref{eq:FALT}), this implies
\beas
\prod_{i=1}^k\lambda_i(f(A)\;f(B))
&\le& \prod_{i=1}^k f_1^{2p}(\lambda_i^{1/2}(AB))\;f_2^{2(1-p)}(\lambda_i^{1/2}(AB)) \\
&=& \prod_{i=1}^k f^{2}(\lambda_i^{1/2}(AB)),
\eeas
i.e.\ $f$ satisfies inequality (\ref{eq:FALT}) as well.
\qed

We have already proven that membership of this class is a
necessary condition for inequality (\ref{eq:FALT}) to hold.
The set of associated functions $F$ for functions in class $\cA$
is the set of concave functions $F$ that satisfy $0\le F'(y)\le 1$ for all $y$ in the domain of $F$ where $F$ is differentiable.
This set is convex, as can be seen from the fact that, for $f\in\cA$,
$F'$ is non-increasing and the range of $F'$ is $[0,1]$. Hence, $F'$ is a convex combination
of step functions $\Phi(b-y)$ (with $\Phi$ the Heaviside step function) and the constant functions $0$ and $1$:
$F'(y)=r+s\int_{(-\infty,+\infty)} \Phi(b-y)d\mu(b)$,
where $r,s\ge0$, $r+s\le1$, and $d\mu$ is a probability measure (normalised positive measure).
Hence, such $F$ have the integral representation 
\be
F(y)=\alpha+ry+s\int_{(-\infty,+\infty)} \min(y,t)\;d\mu(t).
\ee
The additive constant $\alpha$ corresponds to multiplication of $f$ by $e^\alpha$, so we may assume that $\alpha=0$.
Since $r+s\le 1$ it then follows that $f$ is in the geometric convex closure of $f(x)=1$, $f(x)=x$ and 
$f(x)=\min(x,c)$ for $c\in I$ ($c=e^t$).

The next step of the proof is to show that inequality (\ref{eq:FALT}) holds for these extremal functions.
For the functions $f(x)=1$ and $f(x)=x$ this is of course trivial to prove.
Hence let us consider the remaining function $f(x)=\min(x,c)$, with $c\in I$.
As the constant $c$ can be absorbed in the matrices $A$ and $B$, we only need to check the function $f(x)=\min(x,1)$.
The action of this function on a matrix $A$ is to replace any eigenvalue of $A$ that is bigger than $1$ by the value $1$.
I denote this matrix function by $\min(A,1)$.
For this function a stronger inequality can be proven than what is actually needed.
\begin{lemma}
For $A,B\ge0$, and for any $i$
$$
\lambda_i(\min(A,1)\min(B,1))\le \min(\lambda_i(AB),1).
$$
\end{lemma}
\textit{Proof.}
Let $A_1=\min(A,1)$ and $B_1=\min(B,1)$.
We have $A_1\le A$ and $B_1\le B$, so that, using Weyl monotonicity of the eigenvalues twice,
\beas
\lambda_i(A_1B_1)
&=&\lambda_i(A_1^{1/2}B_1A_1^{1/2})\\
&\le& \lambda_i(A_1^{1/2}BA_1^{1/2})\\
&=& \lambda_i(B^{1/2}A_1B^{1/2})\\
&\le& \lambda_i(B^{1/2}AB^{1/2})\\
&=& \lambda_i(AB).
\eeas
This implies also that
$\min(\lambda_i(A_1B_1),1) \le \min(\lambda_i(AB),1)$.

We also have $A_1\le\id$ and $B_1\le\id$, hence by Lemma \ref{lem:a1}
$\min(\lambda_i(A_1B_1),1) = \lambda_i(A_1B_1)$.
\qed

Since the inequality of this lemma implies the weaker log-submajorisation inequality
$$
\lambda(\min(A,1)\min(B,1))\prec_{w(\log)}\min(\lambda(AB),1),
$$
all extremal points of the
class $\cA$ satisfy the inequality (\ref{eq:FALT}).

Finally, by Proposition \ref{prop:convex}
this implies that (\ref{eq:FALT}) holds for all functions in $\cA$, hence
membership of $\cA$ is a sufficient condition. This ends the proof of Theorem \ref{th:FALT}.
\qed

\subsection{Proof of sufficiency for Theorem \ref{cor:RALT}.}
Let $A=g(X)$ and $B=g(Y)$, with $g=f^{-1}$.
Thus, $X=f(A)$ and $Y=f(B)$.
Since $f$ is in $\cA$, $g$ is in $\cB$.
Inequality (\ref{eq:FALT}) then gives
\be
\lambda(XY)\prec_{w(\log)} f^2\left(\sqrt{\lambda(g(X)g(Y))}\right).\label{eq:w}
\ee
The right-hand side features the function $w(x)=f^2(\sqrt{x})$.
Because $f$ is geometrically concave, so is $w$.
The inverse function $w^{-1}$ is given by $w^{-1}(y)=g^2(\sqrt{y})$.
Therefore, $w^{-1}$ is geometrically convex.
Furthermore, because $f'$ is non-negative, $w^{-1}$ is monotonously increasing.

A monotonous convex function preserves the weak majorisation relation (\cite{bhatia}, Corollary II.3.4).
Thus, a monotonous geometrically convex function preserves the log-submajorisation relation.
Hence, when $w^{-1}$ is applied to both sides of (\ref{eq:w}) one obtains
$$
w^{-1}(\lambda(XY))\prec_{w(\log)} \lambda(g(X)g(Y)),
$$
which is (\ref{eq:RALT}).
\qed

\section{Application\label{sec3}}
An interesting application concerns the function $f(x) = 1-\exp(-x)$, which is in class $\cA$.
A simple application of Theorem \ref{th:FALT} leads to an inequality that is complementary to
the famous Golden-Thompson inequality $\trace\exp(A+B)\le\trace\exp(A)\exp(B)$ (where $A$ and $B$ are Hermitian).

In \cite{hiaipetz} (see also \cite{andohiai}) the inequality
$$
\trace(\exp(pA)\#\exp(pB))^{2/p} \le \trace\exp(A+B)
$$
was proven, for every $p>0$ and Hermitian $A$ and $B$.
This is complementary to the Golden-Thompson inequality because it provides a lower bound on $\trace\exp(A+B)$.
The bound obtained below is complementary in a different sense,
as it provides an upper bound on $\trace e^{-A}e^{-B}$ (for \textit{positive} $A$ and $B$).

\begin{theorem}
For $A,B\ge0$, and with
$C=(A^{1/2}BA^{1/2})^{1/2}$,
\be
\trace(e^{-A}e^{-B})\le \trace(e^{-A}+e^{-B}) + \trace(e^{-2C}-2e^{-C}).
\ee
\end{theorem}
\textit{Proof.}
The inequality can be rewritten as
$\trace f(A)f(B)\le \trace f^2(C)$,
with $f(x)=1-e^{-x}$. By Lemma 5 in \cite{ka} $f(x)$ is geometrically concave.
Moreover, $f$ is in $\cA$:
obviously, $f'\ge0$; secondly, $f'(x)=\exp(-x)\le 1/(1+x)$, so that $xf'(x)=x\exp(-x)\le 1-\exp(-x)=f(x)$.

Hence $f$ satisfies the conditions of Theorem \ref{th:FALT}.
The inequality follows immediately from that theorem, as
log-submajorisation implies weak majorisation, and majorisation of the trace, in particular.
\qed

\begin{ack}
Many thanks to F. Hiai for detailed comments on an earlier version of the manuscript.
I am grateful for the hospitality of the University of Ulm, where this paper was completed.
\end{ack}

\end{document}